\documentclass[twoside,11pt,leqno]{article}
\usepackage{amsmath, amssymb, amsthm}

\begin{document}

\frenchspacing

\numberwithin{equation}{section}
\newtheorem{thm}{Theorem}[section]
\newtheorem{lem}[thm]{Lemma}
\newtheorem{prop}[thm]{Proposition}
\newtheorem{cor}[thm]{Corollary}

\theoremstyle{definition}
\newtheorem{defn}{Definition}[section]
\newtheorem{rem}{Remark}[section]


\newcommand{\R}{\mathbb{R}}
\newcommand{\X}{\mathfrak{X}}
\newcommand{\W}{\mathcal{W}}
\newcommand{\s}{\mathfrak{S}}
\newcommand{\g}{\mathfrak{g}}
\newcommand{\Span}{{\rm span}}
\newcommand{\nN}{\norm{N}}
\def\co{\colon\thinspace}
\newcommand{\norm}[1]{\left\Vert#1\right\Vert ^2}
\newcommand{\nJ}{\norm{\nabla J}}
\newcommand{\ad}{{\rm ad}}

\newcommand{\thmref}[1]{Theorem~\ref{#1}}
\newcommand{\dfnref}[1]{Definition~\ref{#1}}


\title{On 4n-dimensional Lie groups as
quasi-K\"ahler manifolds with Killing Norden metric }

\author{Dimitar Mekerov$^1$, Mancho Manev$^2$\footnote{Correspondence author}\\
$^{1,2}$University of Plovdiv, \\
Faculty of Mathematics and Informatics\\
236 Bulgaria Blvd., 4003 Plovdiv, Bulgaria\\
$^{1}$e-mail: mircho@uni-plovdiv.bg \\
$^{2}$e-mail: mmanev@uni-plovdiv.bg  \\
$^{2}$url: http://www.fmi-plovdiv.org/manev/}

\maketitle

\begin{abstract}
A $4n$-parametric family of $4n$-dimensional quasi-K\"ahler
manifolds with Killing Norden metric is constructed on a Lie
group. This family is characterized geometrically.


\end{abstract}

\section{Introduction}

It is a fundamental fact that on an almost complex manifold with
Hermitian metric (almost Hermitian manifold), the action of the
almost complex structure on the tangent space at each point of the
manifold is isometry. There is another type of metric, called a
Norden metric or a $B$-metric on an almost complex manifold, such
as the action of the almost complex structure is anti-isometry
with respect to the metric. Such a manifold is called an almost
complex manifold with Norden metric \cite{GaBo} or with $B$-metric
\cite{GaGrMi}. See also \cite{GrMeDj} for generalized
$B$-manifolds. It is known \cite{GaBo} that these manifolds are
classified into three basic classes $\W_i$ $(i=1,2,3)$ which give
rise eight classes in all.

Among the basic three classes of this classification, the almost
complex structure is nonintegrable only in the class $\W_3$. This
is the class of the so-called \emph{quasi-K\"ahler manifolds with
Norden metric}, which we call briefly \emph{$\W_3$-manifolds}. We
was studied the geometry of manifolds belonging to this class in
\cite{GrMaMe-2}, \cite{MaGrMe-4}, \cite{ManMek-5},
\cite{MekMan-1}, \cite{MeMaGr-3}.

The purpose of the present paper is to exhibit, by construction,
almost complex structures with Norden metric on Lie groups as
$4n$-manifolds, which are of the class $\W_3$.  This
$4n$-parametric family of manifolds is characterized
geometrically.

The case of the initial dimension 4 is considered in
\cite{GrMaMe-2} and \cite{MeMaGr-3}.

\section{Almost complex manifolds with Norden metric}


Let $(M,J,g)$ be a $2n$-dimensional almost complex manifold with
Norden metric, i.~e. $J$ is an almost complex structure and $g$ is
a metric on $M$ such that
\begin{equation}
J^2X=-X, \qquad g(JX,JY)=-g(X,Y)
\end{equation}
for all differentiable vector fields $X$, $Y$ on $M$, i.~e. $X, Y
\in \X(M)$.

The associated metric $\tilde{g}$ of $g$ on $M$ given by
$\tilde{g}(X,Y)=g(X,JY)$ for all $X, Y \in \X(M)$ is a Norden
metric, too. Both metrics are necessarily of signature $(n,n)$.
The manifold $(M,J,\tilde{g})$ is an almost complex manifold with
Norden metric, too.

Further, $X$, $Y$, $Z$, $U$ ($x$, $y$, $z$, $u$, respectively)
will stand for arbitrary differentiable vector fields on $M$
(vectors in $T_pM$, $p\in M$, respectively).

The Levi-Civita connection of $g$ is denoted by $\nabla$. The
tensor filed $F$ of type $(0,3)$ on $M$ is defined by
\begin{equation}\label{F}
F(X,Y,Z)=g\bigl( \left( \nabla_X J \right)Y,Z\bigr).
\end{equation}
It has the following symmetries
\begin{equation}\label{F-prop}
F(X,Y,Z)=F(X,Z,Y)=F(X,JY,JZ).
\end{equation}

Further, let $\{e_i\}$ ($i=1,2,\dots,2n$) be an arbitrary basis of
$T_pM$ at a point $p$ of $M$. The components of the inverse matrix
of $g$ are denoted by $g^{ij}$ with respect to the basis
$\{e_i\}$.

The eight classes of almost complex manifolds with Norden metric
are determined in \cite{GaBo} according to the properties of $F$.
The three basic classes are given as it follows:
\[
\begin{array}{l}
\W_1\co F(x,y,z)=\frac{1}{4n} \left\{
g(x,y)\theta(z)+g(x,z)\theta(y)\right. \\[4pt]
\phantom{\mathcal{W}_1\co F(x,y,z)=\frac{1}{4n} }\left.
    +g(x,J y)\theta(J z)
    +g(x,J z)\theta(J y)\right\};\\[4pt]
\W_2\co \mathop{\s} \limits_{x,y,z}
F(x,y,J z)=0,\quad \theta=0;\\[8pt]
\W_3\co \mathop{\s} \limits_{x,y,z} F(x,y,z)=0,
\end{array}
\]
where $\s $ is the cyclic sum over three arguments and
$\theta(z)=g^{ij}F(e_i,e_j,z)$.
The special class $\W_0$ of the K\"ahler manifolds with Norden
metric belonging to any other class is determined by the condition
$F=0$.


The curvature tensor field $R$ of $\nabla$ is $R(X,Y)Z=\nabla_X
\nabla_Y Z - \nabla_Y \nabla_X Z -
    \nabla_{[X,Y]}Z$ and the corresponding tensor field of type $(0,4)$ is determined
    by $R(X,Y,Z,U)=g(R(X,Y)Z,U)$. The Ricci tensor $\rho$ and the
scalar curvature $\tau$ are defined as usual by
\begin{equation}\label{rho-tau}
    \rho(y,z)=g^{ij}R(e_i,y,z,e_j),\qquad \tau=g^{ij}\rho(e_i,e_j).
\end{equation}

It is well-known that the Weyl tensor $W$ on a $m$-dimensional
pseudo-Riemannian manifold ($m\geq 3$) is given by
\begin{equation}\label{W}
    W=R-\frac{1}{m-2}\Bigl(\psi_1(\rho)-\frac{\tau}{m-1}\pi_1\Bigr),
\end{equation}
where
\[
    \begin{array}{l}
      \psi_1(\rho)(x,y,z,u)=g(y,z)\rho(x,u)-g(x,z)\rho(y,u) \\
      \phantom{\psi_1(\rho)(x,y,z,u)}+\rho(y,z)g(x,u)-\rho(x,z)g(y,u); \\
      \pi_1=\frac{1}{2}\psi_1(g)=g(y,z)g(x,u)-g(x,z)g(y,u). \\
    \end{array}
\]

Moreover, for $m\geq 4$ the Weyl tensor $W$ is zero if and only if
the manifold is conformally flat.

Let $\alpha=\{x,y\}$ be a non-degenerate 2-plane spanned by
vectors $x, y \in T_pM$, $p\in M$. It means that
$\pi_1(x,y,y,x)=g(x,x)g(y,y)-g(x,y)^2 \neq 0$). Then, it is known,
the sectional curvature of $\alpha$ is defined by the following
equation
\begin{equation}\label{k}
    k(\alpha)=k(x,y)=\frac{R(x,y,y,x)}{\pi_1(x,y,y,x)}.
\end{equation}

The basic sectional curvatures in $T_pM$ with an almost complex
structure and a Norden metric $g$ are:
\begin{itemize}
    \item \emph{holomorphic sectional curvatures} if $J\alpha=\alpha$;
    \item \emph{totally real sectional curvatures} if
    $J\alpha\perp\alpha$ with respect to $g$.
\end{itemize}

In \cite{GrDjMe}, a \emph{holomorphic bisectional curvature}
$h(x,y)$ for a pair of holomorphic 2-planes $\alpha_1=\{x,Jx\}$
and $\alpha_2=\{y,Jy\}$ is defined by
\begin{equation}\label{h}
    h(x,y)=-\frac{R(x,Jx,y,Jy)}
    {\sqrt{\pi_1(x,Jx,x,Jx)\pi_1(y,Jy,y,Jy)}},
\end{equation}
where $x$, $y$ do not lie along the totally isotropic directions,
i.~e. both of the couples $\bigl(g(x,x), g(x,Jx)\bigr)$ and
$\bigl(g(y,y), g(y,Jy)\bigr)$ are different from the couple
$\left(0,0\right)$. The holomorphic bisectional curvature is
invariant with respect to the basis of the 2-planes $\alpha_1$ and
$\alpha_2$. In particular, if $\alpha_1=\alpha_2$, then the
holomorphic bisectional curvature coincides with the holomorphic
sectional curvature of the 2-plane $\alpha_1=\alpha_2$.


 The square norm $\nJ$ of $\nabla J$ is defined in
\cite{GRMa} by
\[
    \nJ=g^{ij}g^{kl}
    g\bigl(\left(\nabla_{e_i} J\right)e_k,\left(\nabla_{e_j}
    J\right)e_l\bigr).
\]

Having in mind the definition \eqref{F} of the tensor $F$ and the
properties \eqref{F-prop}, we obtain the following equation for
the square norm of $\nabla J$
\begin{equation}\label{snormF}
    \nJ=g^{ij}g^{kl}g^{pq}F_{ikp}F_{jlq},
\end{equation}
where $F_{ikp}=F(e_i,e_k,e_p)$.

An almost complex manifold with Norden metric satisfying the
condition $\nJ=0$ is called an \emph{isotropic K\"ahler manifold
with Norden metric} \cite{MekMan-1}. It is clear, if a manifold
belongs to the class $\W_0$, then it is isotropic K\"ahlerian but
the inverse statement is not always true.


\section{A Lie group as a $4n$-dimensional $\W_3$-mani\-fold}

Let $V$ be a $4n$-dimensional vector space and let us consider the
structure of the Lie algebra $\g$ defined by the brackets
$[E_i,E_j]=C_{ij}^kE_k$, where $\{E_1,E_2,\dots,E_{4n}\}$ is a
basis of $V$ and $C_{ij}^k\in \R$.

Let $G$ be the associated connected Lie group and
\(\{X_1,X_2,\dots,X_{4n}\}\) be a global basis of left invariant
vector fields. Then the Jacobi identity holds:
\begin{equation}\label{Jac}
    \mathop{\s} \limits_{X_i,X_j,X_k}
    \bigl[[X_i,X_j],X_k\bigr]=0.
\end{equation}

Next we define an almost complex structure $J$ by the conditions
\begin{equation}\label{J4n}
\begin{array}{ll}
JX_{4\alpha-3}=X_{4\alpha-1},\quad &
JX_{4\alpha-2}=X_{4\alpha},\quad\\[4pt]
JX_{4\alpha-1}=-X_{4\alpha-3},\quad & JX_{4\alpha}=-X_{4\alpha-2},
\end{array}
\end{equation}
where $\alpha\in\{1,2,\dots,n\}$.

Let us consider the left invariant metric $g$ defined by the
following way
\begin{equation}\label{g}
    \begin{split}
  g(X_{4\alpha-3},X_{4\alpha-3})&=g(X_{4\alpha-2},X_{4\alpha-2})=-g(X_{4\alpha-1},X_{4\alpha-1})\\[4pt]
  &=-g(X_{4\alpha},X_{4\alpha})=1,
\\[4pt]
g(X_i,X_j)&=0\quad \text{for}\quad i\neq j. \\
\end{split}
\end{equation}
The introduced metric is a Norden metric because of \eqref{J4n}.

In this way, the induced $4n$-dimensional manifold $(G,J,g)$ is an
almost complex manifold with Norden metric, in short \emph{almost
Norden manifold}.

From this point on, until the end of this paper we shall consider
almost Norden manifolds $(G,J,g)$ with \emph{Killing metric} $g$.
This means that $g$ satisfies the following condition for
arbitrary $X,Y,Z \in \g$
\[
    g\left([X,Y],Z\right)+g\left([X,Z],Y\right)=0.
\]

In \cite{MaGrMe-4} it is known that each almost Norden manifold
with Killing metric is a locally symmetric $\W_3$-manifold.
Moreover, the following formulae are valid:
\begin{equation}\label{nabla[]}
    \nabla_{X_i} X_j=\frac{1}{2}[X_i,X_j],
\end{equation}
\begin{equation}\label{F[]}
    F(X_i,X_j,X_k)=\frac{1}{2}\Bigl\{g\bigl( [X_i, JX_j],X_k\bigr) -
    g\bigl([X_i, X_j],JX_k \bigr) \Bigr\},
\end{equation}
\begin{equation}\label{R[]}
   R(X_i,X_j,X_k,X_l)=-\frac{1}{4}g\Bigl([X_i,X_j],[X_k,X_l]\Bigr),
\end{equation}
\begin{equation}\label{ort}
    [X_i,X_j]\bot\;\Span\{X_i,X_j\},
\end{equation}
where $i,j,k,l\in\{1,2,\dots,4n\}$.

Since $g$ is a Killing metric, the structural constants $C_{ij}^k$
are specialized such that the commutators have the following
decompositions:
\begin{equation}\label{[]4}
\begin{array}{l}
[X_{4\alpha-3},X_{4\alpha-1}]=\lambda_{4\alpha-2} X_{4\alpha-2}+\lambda_{4\alpha} X_{4\alpha},\\[4pt]
[X_{4\alpha-2},X_{4\alpha}]=\lambda_{4\alpha-3} X_{4\alpha-3}+\lambda_{4\alpha-1} X_{4\alpha-1},\\[4pt]
[X_{4\alpha-2},X_{4\alpha-1}]=-\lambda_{4\alpha-2} X_{4\alpha-3}-\lambda_{4\alpha-1} X_{4\alpha},\\[4pt]
[X_{4\alpha-1},X_{4\alpha}]=-\lambda_{4\alpha} X_{4\alpha-3}+\lambda_{4\alpha-1} X_{4\alpha-2},\\[4pt]
[X_{4\alpha},X_{4\alpha-3}]=\lambda_{4\alpha-3} X_{4\alpha-2}+\lambda_{4\alpha} X_{4\alpha-1},\\[4pt]
[X_{4\alpha-2},X_{4\alpha-3}]=-\lambda_{4\alpha-2} X_{4\alpha-1}+\lambda_{4\alpha-3} X_{4\alpha},\\[4pt]
\end{array}
\end{equation}
where $\lambda_{4\alpha-i}\in\R$ $(i = 1, 2, 3, 4; \alpha = 1, 2,
\dots, n)$. The other commutators are zero.

By direct verification we prove that the commutators from
\eqref{[]4} satisfy the Jacobi identity \eqref{Jac}. The Lie
groups $G$ thus obtained are a family which is characterized by
$4n$ real parameters $\lambda_{4\alpha-i}$.

Vice versa, let the condition \eqref{[]4} is valid for an almost
Norden manifold $(G,J,g)$ with a structure $J$ and a metric $g$
determined by \eqref{J4n} and  \eqref{g}, respectively. Then we
verify directly that $g$ is a Killing metric, i.e. $(G,J,g)$ is
locally symmetric $\W_3$-manifold.

Therefore we establish the truthfulness of the following
\begin{thm}\label{Thm-4n}
Let $(G,J,g)$ be a $4n$-dimensional almost Norden manifold, where
$G$ is a connected Lie group with corresponding Lie algebra $\g$
determined by the global basis of left invariant vector fields
$\{X_1,X_2,\dots,X_{4n}\}$; $J$ is an almost complex structure
defined by \eqref{J4n} and $g$ is a Norden metric determined by
\eqref{g}. Then $(G,J,g)$ is a $\W_3$-manifold with Killing metric
$g$ if and only if $G$ belongs to the $4n$-parametric family of
Lie groups determined by the conditions \eqref{[]4}.
\end{thm}

\section{Geometric characteristics of the constructed
ma\-ni\-fold}

Let $(G,J,g)$ be the $4n$-dimensional quasi-K\"ahler manifold with
Norden metric introduced in the previous section. Let us introduce
the following index denotations: $\bar{1}=4\alpha-3$,
$\bar{2}=4\alpha-2$, $\bar{3}=4\alpha-1$, $\bar{4}=4\alpha$ for
any fixed $\alpha\in\{1,2,\dots,n\}$.


Having in mind \eqref{F[]}, \eqref{J4n}, \eqref{g} and
\eqref{[]4}, we obtain immediately the nonzero components of the
tensor $F$ as follows:

\begin{equation}\label{Fijk}
\begin{array}{r}
-F_{\bar{1}\bar{2}\bar{2}}=-F_{\bar{1}\bar{4}\bar{4}}=2F_{\bar{2}\bar{1}\bar{2}}
=2F_{\bar{2}\bar{3}\bar{4}}
=2F_{\bar{4}\bar{1}\bar{4}}=-2F_{\bar{4}\bar{2}\bar{3}}
=\lambda_{\bar{1}},\\[4pt]
2F_{\bar{1}\bar{1}\bar{2}}=2F_{\bar{1}\bar{3}\bar{4}}
=-2F_{\bar{2}\bar{1}\bar{1}}=-2F_{\bar{2}\bar{3}\bar{3}}=-2F_{\bar{3}\bar{1}\bar{4}}=2F_{\bar{3}\bar{2}\bar{3}}
=\lambda_{\bar{2}},\\[4pt]
2F_{\bar{2}\bar{1}\bar{4}}=-2F_{\bar{2}\bar{2}\bar{3}}
=F_{\bar{3}\bar{2}\bar{2}}=F_{\bar{3}\bar{4}\bar{4}}=-2F_{\bar{4}\bar{1}\bar{2}}
=-2F_{\bar{4}\bar{3}\bar{4}}=\lambda_{\bar{3}},\\[4pt]
-2F_{\bar{1}\bar{1}\bar{4}}=2F_{\bar{1}\bar{2}\bar{3}}
=-2F_{\bar{3}\bar{1}\bar{2}}=-2F_{\bar{3}\bar{3}\bar{4}}
=F_{\bar{4}\bar{1}\bar{1}}=F_{\bar{4}\bar{3}\bar{3}}=\lambda_{\bar{4}}.
\end{array}
\end{equation}
The other nonzero components are obtained from the property
$F_{\bar{i}\bar{j}\bar{k}}=F_{\bar{i}\bar{k}\bar{j}}$.


 Let $N$ be the Nijenhuis tensor of the almost complex
structure $J$ on $G$, i.e.
\[
    N(X,Y)=[X,Y]+J[JX,Y]+J[X,JY]-[JX,JY], \quad X, Y \in \g.
\]
Having in mind \eqref{[]4} and \eqref{J4n} we obtain the nonzero
components $N_{\bar{i}\bar{j}}=N(X_{\bar{i}},X_{\bar{j}})$ as
follows
\begin{equation}\label{Nij}
\begin{array}{l}
    N_{\bar{1}\bar{2}}=-N_{\bar{3}\bar{4}}
    =2\left(\lambda_{\bar{4}} X_{\bar{1}}-\lambda_{\bar{3}}X_{\bar{2}}+\lambda_{\bar{2}}
    X_{\bar{3}}-\lambda_{\bar{1}}X_{\bar{4}}\right),\\[4pt]
    N_{\bar{1}\bar{4}}=-N_{\bar{2}\bar{3}}=2\left(\lambda_{\bar{2}} X_{\bar{1}}-\lambda_{\bar{1}}
    X_{\bar{2}}-\lambda_{\bar{4}} X_{\bar{3}}+\lambda_{\bar{3}} X_{\bar{4}} \right).\\[4pt]
\end{array}
\end{equation}
The other nonzero components are obtained from the property
$N_{\bar{i}\bar{j}}=-N_{\bar{j}\bar{i}}$.
Hence its square norm $\nN=g^{ik}g^{ks}g(N_{ij},N_{ks})$ for
$i,j,k,s\in\{1,2,\dots,4n\}$ has the form
\begin{equation}\label{nN}
    \nN=-32\sum_{\alpha=1}^n \left(\lambda_{4\alpha-3}^2+\lambda_{4\alpha-2}^2
    -\lambda_{4\alpha-1}^2-\lambda_{4\alpha}^2\right),
\end{equation}
where the inverse matrix of $g$ has the form
\begin{equation}\label{g^ij}
    \left(g^{ij}\right)=\left(%
\begin{array}{cccc}
  \tilde{E}_4 & 0 & \dots & 0 \\
  0 & \tilde{E}_4 & \dots & 0 \\
  \dots & \dots & \dots & \dots \\
  0 & 0 & 0 & \tilde{E}_4 \\
\end{array}%
 \right), \quad
    \tilde{E}_4=\left(%
\begin{array}{cccc}
  1 & 0 & 0 & 0 \\
  0 & 1 & 0 & 0 \\
  0 & 0 & -1 & 0 \\
  0 & 0 & 0 & -1 \\
\end{array}%
\right)
\end{equation}


According to \eqref{J4n}, \eqref{g}, \eqref{[]4}  and
\eqref{g^ij}, from \eqref{snormF} we obtain the square norm of
$\nabla J$ as
\begin{equation}\label{snorm4}
    \nJ=4\sum_{\alpha=1}^n \left(\lambda_{4\alpha-3}^2+\lambda_{4\alpha-2}^2
    -\lambda_{4\alpha-1}^2-\lambda_{4\alpha}^2\right).
\end{equation}


From \eqref{R[]} and \eqref{[]4} we get the nonzero components of
$R$ as follows
\begin{equation}\label{Rijks}
\begin{array}{ll}
    R_{{\bar{1}}{\bar{2}}{\bar{2}}{\bar{1}}}=-\frac{1}{4}\left(\lambda_{\bar{1}}^2+\lambda_{\bar{2}}^2\right),\quad
    &
    R_{{\bar{1}}{\bar{3}}{\bar{3}}{\bar{1}}}=\frac{1}{4}\left(\lambda_{\bar{2}}^2-\lambda_{\bar{4}}^2\right),\\[4pt]
    R_{{\bar{1}}{\bar{4}}{\bar{4}}{\bar{1}}}=-\frac{1}{4}\left(\lambda_{\bar{1}}^2-\lambda_{\bar{4}}^2\right),\quad
    &
    R_{{\bar{2}}{\bar{3}}{\bar{3}}{\bar{2}}}=\frac{1}{4}\left(\lambda_{\bar{2}}^2-\lambda_{\bar{3}}^2\right),\\[4pt]
    R_{{\bar{2}}{\bar{4}}{\bar{4}}{\bar{2}}}=\frac{1}{4}\left(\lambda_{\bar{1}}^2-\lambda_{\bar{3}}^2\right),\quad
    &
    R_{{\bar{3}}{\bar{4}}{\bar{4}}{\bar{3}}}=\frac{1}{4}\left(\lambda_{\bar{3}}^2+\lambda_{\bar{4}}^2\right),\\[4pt]
    R_{{\bar{1}}{\bar{3}}{\bar{4}}{\bar{1}}}=R_{{\bar{2}}{\bar{3}}{\bar{4}}{\bar{2}}}=-\frac{1}{4}\lambda_{\bar{1}}\lambda_{\bar{2}},\quad
    &
    R_{{\bar{2}}{\bar{1}}{\bar{3}}{\bar{2}}}=-R_{{\bar{4}}{\bar{1}}{\bar{3}}{\bar{4}}}=\frac{1}{4}\lambda_{\bar{1}}\lambda_{\bar{3}},\\[4pt]
    R_{{\bar{1}}{\bar{2}}{\bar{3}}{\bar{1}}}=-R_{{\bar{4}}{\bar{2}}{\bar{3}}{\bar{4}}}=\frac{1}{4}\lambda_{\bar{1}}\lambda_{\bar{4}},\quad
    &
    R_{{\bar{2}}{\bar{1}}{\bar{4}}{\bar{2}}}=-R_{{\bar{3}}{\bar{1}}{\bar{4}}{\bar{3}}}=\frac{1}{4}\lambda_{\bar{2}}\lambda_{\bar{3}},\\[4pt]
    R_{{\bar{1}}{\bar{2}}{\bar{4}}{\bar{1}}}=-R_{{\bar{3}}{\bar{2}}{\bar{4}}{\bar{3}}}=\frac{1}{4}\lambda_{\bar{2}}\lambda_{\bar{4}},\quad
    &
    R_{{\bar{3}}{\bar{1}}{\bar{2}}{\bar{3}}}=R_{{\bar{4}}{\bar{1}}{\bar{2}}{\bar{4}}}=\frac{1}{4}\lambda_{\bar{3}}\lambda_{\bar{4}}.\\[4pt]
\end{array}
\end{equation}
The other nonzero components of $R$ are get from the properties
$R_{\bar{i}\bar{j}\bar{k}\bar{s}}=R_{\bar{k}\bar{s}\bar{i}\bar{j}}$
and
$R_{\bar{i}\bar{j}\bar{k}\bar{s}}=-R_{\bar{j}\bar{i}\bar{k}\bar{s}}=-R_{\bar{i}\bar{j}\bar{s}\bar{k}}$.


Having in mind \eqref{rho-tau}, \eqref{g^ij} and \eqref{Rijks}, we
obtain the components
$\rho_{\bar{i}\bar{j}}=\rho(X_{\bar{i}},X_{\bar{j}})$ of the Ricci
tensor $\rho$ and the scalar curvature $\tau$ as follows:
\begin{equation}\label{rho_ij}
\begin{array}{c}
\begin{array}{ll}
    \rho_{{\bar{1}}{\bar{1}}}=-\frac{1}{2}\left(\lambda_{\bar{1}}^2+\lambda_{\bar{2}}^2-\lambda_{\bar{4}}^2\right),\quad
    &
    \rho_{{\bar{2}}{\bar{2}}}=-\frac{1}{2}\left(\lambda_{\bar{1}}^2+\lambda_{\bar{2}}^2-\lambda_{\bar{3}}^2\right),\\[4pt]
    \rho_{{\bar{3}}{\bar{3}}}=\frac{1}{2}\left(\lambda_{\bar{2}}^2-\lambda_{\bar{3}}^2-\lambda_{\bar{4}}^2\right),\quad
    &
    \rho_{{\bar{4}}{\bar{4}}}=\frac{1}{2}\left(\lambda_{\bar{1}}^2-\lambda_{\bar{3}}^2-\lambda_{\bar{4}}^2\right),\\[4pt]
    \rho_{{\bar{1}}{\bar{2}}}=\rho_{{\bar{2}}{\bar{1}}}=-\frac{1}{2}\lambda_{\bar{3}}\lambda_{\bar{4}},\quad
    &
    \rho_{{\bar{1}}{\bar{3}}}=\rho_{{\bar{3}}{\bar{1}}}=\frac{1}{2}\lambda_{\bar{1}}\lambda_{\bar{3}},\quad\\[4pt]
    \rho_{{\bar{1}}{\bar{4}}}=\rho_{{\bar{4}}{\bar{1}}}=\frac{1}{2}\lambda_{\bar{2}}\lambda_{\bar{3}},\quad
    &
    \rho_{{\bar{2}}{\bar{3}}}=\rho_{{\bar{3}}{\bar{2}}}=\frac{1}{2}\lambda_{\bar{1}}\lambda_{\bar{4}},\quad\\[4pt]
    \rho_{{\bar{2}}{\bar{4}}}=\rho_{{\bar{4}}{\bar{2}}}=\frac{1}{2}\lambda_{\bar{2}}\lambda_{\bar{4}},\quad
    &
    \rho_{{\bar{3}}{\bar{4}}}=\rho_{{\bar{4}}{\bar{3}}}=-\frac{1}{2}\lambda_{\bar{1}}\lambda_{\bar{2}};\\[4pt]
\end{array}
\end{array}
\end{equation}
\begin{equation}\label{tau}
    \tau=-\frac{3}{2}\sum_{\alpha=1}^n \left(\lambda_{4\alpha-3}^2+\lambda_{4\alpha-2}^2
    -\lambda_{4\alpha-1}^2-\lambda_{4\alpha}^2\right).
\end{equation}


Taking into account \eqref{g}, \eqref{Rijks}, \eqref{rho_ij},
\eqref{tau} and \eqref{W} for $m=4n$, we establish that the Weyl
tensor vanishes. Then $(G,J,g)$ is a conformally flat manifold.


For the sectional curvatures
$k_{\bar{i}\bar{j}}=k(\alpha_{\bar{i}\bar{j}})$ of the basic
2-planes
$\alpha_{\alpha_{\bar{i}\bar{j}}}=\{X_{\bar{i}},X_{\bar{j}}\}$,
according to \eqref{k}, \eqref{g} and \eqref{Rijks}, we have:
\begin{equation}\label{k_ik}
    \begin{array}{ll}
    k_{\bar{1}\bar{3}}=-\frac{1}{4}\left(\lambda_{\bar{2}}^2-\lambda_{\bar{4}}^2\right),\quad
    &
    k_{\bar{2}\bar{4}}=-\frac{1}{4}\left(\lambda_{\bar{1}}^2-\lambda_{\bar{3}}^2\right),\\[4pt]
    k_{\bar{1}\bar{2}}=-\frac{1}{4}\left(\lambda_{\bar{1}}^2+\lambda_{\bar{2}}^2\right),\quad
    &
    k_{\bar{1}\bar{4}}=-\frac{1}{4}\left(\lambda_{\bar{1}}^2-\lambda_{\bar{4}}^2\right),\\[4pt]
    k_{\bar{2}\bar{3}}=-\frac{1}{4}\left(\lambda_{\bar{2}}^2-\lambda_{\bar{3}}^2\right),\quad
    &
    k_{\bar{3}\bar{4}}=\frac{1}{4}\left(\lambda_{\bar{3}}^2+\lambda_{\bar{4}}^2\right).\\[4pt]
\end{array}
\end{equation}


The obtained geometric characteristics of the considered manifold
we generalize in the following
\begin{thm}\label{thm1-end}
Let $(G,J,g)$ be a $4n$-dimensional almost Norden manifold, where
$G$ is a connected Lie group with corresponding Lie algebra $\g$
determined by the global basis of left invariant vector fields
$\{X_1,X_2,\dots,X_{4n}\}$; $J$ is an almost complex structure
defined by \eqref{J4n} and $g$ is a Norden metric determined by
\eqref{g}. Then
   \begin{enumerate}
    \renewcommand{\labelenumi}{(\roman{enumi})}
    \item
    $(G,J,g)$ is a locally symmetric conformally flat $\W_3$-manifold with Kill\-ing metric $g$;
    \item
    The nonzero components of the basic tensor $F$, the Nijenhuis tensor $N$, the curvature tensor $R$
    and the Ricci tensor $\rho$ are
    \eqref{Fijk}, \eqref{Nij}, \eqref{Rijks} and \eqref{rho_ij}, respectively;
    \item
    The square norms of the Nijenhuis tensor $N$ and $\nabla J$ are \eqref{nN} and \eqref{snorm4}, respectively;
    \item
    The scalar curvature $\tau$ and the sectional curvatures $k_{\bar{i}\bar{j}}$
    of the basic 2-planes are \eqref{tau} and \eqref{k_ik}, respectively.
    \end{enumerate}
\end{thm}

The last theorem implies immediately the following
\begin{cor}\label{thm2-end}
Let $(G,J,g)$ be a $4n$-dimensional almost Norden manifold, where
$G$ is a connected Lie group with corresponding Lie algebra $\g$
determined by the global basis of left invariant vector fields
$\{X_1,X_2,\dots,X_{4n}\}$; $J$ is an almost complex structure
defined by \eqref{J4n} and $g$ is a Norden metric determined by
\eqref{g}. Then the following propositions are equivalent:
\begin{enumerate}
    \renewcommand{\labelenumi}{(\roman{enumi})}
    \item
    $(G,J,g)$ is an isotropic K\"ahler manifold;
    \item
    $(G,J,g)$ is a scalar flat manifold;
    \item
    The Nijenhuis tensor is isotropic;
    \item
    The condition $\sum_{\alpha=1}^n \left(\lambda_{4\alpha-3}^2+\lambda_{4\alpha-2}^2
    -\lambda_{4\alpha-1}^2-\lambda_{4\alpha}^2\right)=0$
    holds.
    \end{enumerate}
\end{cor}

The condition (iv) of the last theorem means that the set of
vectors with the coordinates
$(\lambda_1,\lambda_2,\dots,\lambda_{4n})$ at an arbitrary point
$p\in G$ describes the isotropic cone in $T_pG$ with respect to
the Norden metric $g$.

Let us remark that the 2-planes $\alpha_{\bar{1}\bar{3}}$ and
$\alpha_{\bar{2}\bar{4}}$ are holomorphic 2-planes and the
2-planes $\alpha_{\bar{1}\bar{2}}$, $\alpha_{\bar{1}\bar{4}}$,
$\alpha_{\bar{2}\bar{3}}$, $\alpha_{\bar{3}\bar{4}}$ are totally
real 2-planes.
Taking into account \eqref{h}, \eqref{g} and \eqref{Rijks}, we
obtain that the holomorphic bisectional curvature of the unique
pair of basis holomorphic 2-planes $\{\alpha_{\bar{1}\bar{3}},
\alpha_{\bar{2}\bar{4}}\}$ vanishes.
Moreover, the equalities \eqref{k_ik} imply the following

\begin{thm}\label{thm3-end}
Let $(G,J,g)$ be a $4n$-dimensional almost Norden manifold, where
$G$ is a connected Lie group with corresponding Lie algebra $\g$
determined by the global basis of left invariant vector fields
$\{X_1,X_2,\dots,X_{4n}\}$; $J$ is an almost complex structure
defined by \eqref{J4n} and $g$ is a Norden metric determined by
\eqref{g}. Then
   \begin{enumerate}
   \renewcommand{\labelenumi}{(\roman{enumi})}
    \item
    $(G,J,g)$ is of constant holomorphic sectional curvatures
    iff
    \[\lambda_{\bar{1}}^2+\lambda_{\bar{4}}^2=\lambda_{\bar{2}}^2+\lambda_{\bar{3}}^2;\]
    \item
    $(G,J,g)$ does not admit constant totally real sectional
    curvatures.
    \end{enumerate}
\end{thm}

\end{document}